\documentclass[10pt,reqno]{amsart}
\usepackage[utf8]{inputenc}
\usepackage[english]{babel}
\usepackage{amsmath, amssymb, amsthm, comment}
\usepackage{mathrsfs}
\usepackage{enumitem}
\usepackage{graphicx}
\usepackage{tikz-cd}
\usepackage{xcolor}
\numberwithin{equation}{section} 

\usepackage{diagbox} 
\usepackage{tensor}  

\usepackage{hyperref} 

\usepackage[alphabetic]{amsrefs} 

\theoremstyle{plain}
\newtheorem{thm}{Theorem}[section]
\newtheorem{prop}[thm]{Proposition}
\newtheorem{lemma}[thm]{Lemma}

\theoremstyle{definition}
\newtheorem{defin}[thm]{Definition}

\newtheorem{rmk}[thm]{Remark}

\allowdisplaybreaks

\def\XXint#1#2#3{{\setbox0=\hbox{$#1{#2#3}{\int}$ }
\vcenter{\hbox{$#2#3$ }}\kern-.6\wd0}}

\usepackage{scalerel,stackengine}
\stackMath
\newcommand\hhat[1]{%
\savestack{\tmpbox}{\stretchto{%
  \scaleto{%
    \scalerel*[\widthof{\ensuremath{#1}}]{\kern.1pt\mathchar"0362\kern.1pt}%
    {\rule{0ex}{\textheight}}
  }{\textheight}%
}{2.4ex}}%
\stackon[-6.9pt]{#1}{\tmpbox}%
}

\title[Determining complete Riemannian manifolds from local space-time]{Determining complete Riemannian manifolds\\ from local space-time}

\author{Matti Lassas}

\author{Ruirui Wu}

\begin{document}

\begin{abstract}
We consider a space--time fractional diffusion equation on a complete
Riemannian manifold $(M,g)$ of dimension $n\geq 2$. We prove that,
for a suitably chosen source supported in a relatively compact open
subset $V\Subset M$, a single measurement of the corresponding solution
on $V$ over a finite time interval determines the manifold up to a
Riemannian isometry. The proof combines the source construction with
time analyticity to recover the local source-to-solution map. We then
use the Laplace transform and Stone's formula to recover the local
spectral measure, reducing the geometric determination to the
corresponding inverse problem for the wave equation.
\end{abstract}

\maketitle

\section{Introduction}
We study the recovery of a Riemannian manifold from a single local
measurement of a space--time fractional diffusion process.
Let $(M,g)$ be a connected, complete, smooth Riemannian manifold
without boundary, of dimension $n\geq2$, and let $V\Subset M$
be a nonempty open subset with smooth boundary.
Fix $\alpha,\beta\in(0,1]$ and an observation time $T>0$.

Let $-\Delta_g$ denote the nonnegative self-adjoint realization
of the Laplace--Beltrami operator on $L^2(M)$.
Its fractional power is defined by
\[
    (-\Delta_g)^\alpha
    =
    \int_{[0,\infty)}
        \lambda^\alpha\,dE_g(\lambda),
\]
where $E_g$ is the spectral measure of $-\Delta_g$.
For a Banach space $X$ and a function
$u\in C^1([0,\infty);X)$, the Caputo derivative of order $\beta$
is given by
\[
    \partial_t^\beta u(t)
    :=
    \begin{cases}
    \displaystyle
    \frac{1}{\Gamma(1-\beta)}
    \int_0^t
        (t-\tau)^{-\beta}\,\partial_\tau u(\tau)\,d\tau,
        & 0<\beta<1,\\[1.2em]
    \displaystyle
    \partial_tu(t),
        & \beta=1.
    \end{cases}
\]
Here $\Gamma$ is the Gamma function, and the integral is understood
in the Bochner sense.

Denote by $e_V$ extension by zero from $V$ to $M$, and by $r_V$
restriction from $M$ to $V$.
For
\[
    f\in C_c^2\bigl((0,\infty);L^2(V)\bigr),
\]
we consider the initial value problem
\begin{equation}\label{eq:fractional-diffusion}
    \begin{cases}
        \partial_t^\beta u_f+(-\Delta_g)^\alpha u_f=e_Vf
            & \text{in }M\times(0,\infty),\\
        u_f(\cdot,0)=0
            & \text{in }M.
    \end{cases}
\end{equation}
The well-posedness result in
Section~\ref{subsec:source-to-solution} allows us to define the
local source-to-solution operator
\[
    L_V:
    C_c^2\bigl((0,\infty);L^2(V)\bigr)
    \longrightarrow
    C\bigl([0,\infty);L^2(V)\bigr),
    \qquad
    L_Vf:=r_Vu_f.
\]
Our main result shows that the response to a suitably chosen
source, observed only on $V$ during a finite time interval,
determines the entire manifold.

\begin{thm}\label{thm:single-measurement}
Let $(M,g)$ be a connected, complete, smooth Riemannian manifold
without boundary, of dimension $n\geq2$, and let $V\Subset M$
be a nonempty open subset with smooth boundary.
Fix $\alpha,\beta\in(0,1]$ and $T>0$.
Then there exists a source
\[
    h\in C_c^\infty\bigl((0,T);L^2(V)\bigr)
\]
such that the single measurement
\[
    \left(V,u_h|_{V\times[0,T)}\right)
\]
determines $(M,g)$ up to a Riemannian isometry.
Here $u_h$ is the solution of
\eqref{eq:fractional-diffusion} with source $h$.

More precisely, let $(\widetilde M,\widetilde g)$ be another
connected, complete, smooth Riemannian manifold without boundary,
of dimension $n$, and let
$\widetilde V\Subset\widetilde M$ be a nonempty open subset
with smooth boundary.
Suppose that there exists a diffeomorphism
\[
    \theta:
    \overline{\widetilde V}\longrightarrow\overline V
\]
such that
\[
    \theta^*
    \left(u_h|_{V\times[0,T)}\right)
    =
    \widetilde u_{\theta^*h}
    |_{\widetilde V\times[0,T)}.
\]
The pullback acts only on the spatial variable:
\[
    (\theta^*F)(x,t):=F(\theta(x),t),
\]
and $\widetilde u_{\theta^*h}$ denotes the solution of
\eqref{eq:fractional-diffusion} on
$(\widetilde M,\widetilde g)$ with source $\theta^*h$,
zero initial data, and the same fractional orders
$\alpha$ and $\beta$.
Then there exists a Riemannian isometry
\[
    \Phi:(M,g)\longrightarrow
    (\widetilde M,\widetilde g).
\]
\end{thm}

The proof begins with a continuation property for the response
to a source of compact temporal support.
After the source has ceased, the solution is analytic in time.
Consequently, an observation interval extending beyond the support
of the source determines the response at all later times.

\begin{prop}\label{prop:finite-to-infinite-time}
Let $(M,g)$, $(\widetilde M,\widetilde g)$, $V$,
$\widetilde V$, and $\theta$ satisfy the assumptions of
Theorem~\ref{thm:single-measurement}.
Let $T_0>0$ and
\[
    f\in C_c^2\bigl((0,T_0);L^2(V)\bigr),
\]
and extend $f$ by zero to $(0,\infty)$.
If
\[
    \theta^*
    \left(u_f|_{V\times[0,T_0)}\right)
    =
    \widetilde u_{\theta^*f}
    |_{\widetilde V\times[0,T_0)},
\]
then
\[
    \theta^*
    \left(u_f|_{V\times[0,\infty)}\right)
    =
    \widetilde u_{\theta^*f}
    |_{\widetilde V\times[0,\infty)}.
\]
Equivalently,
\[
    \theta^*(L_Vf)
    =
    \widetilde L_{\widetilde V}(\theta^*f)
    \qquad\text{on }[0,\infty).
\]
\end{prop}

To recover the full source-to-solution operator, the source $h$ is constructed in \cite{HelinLassasYlinenZhang2019}
as a countable superposition of separated temporal pulses whose
spatial profiles form a dense family in $L^2(V)$.
Causality and the continuation property allow us to recover
the responses to the individual pulses from the single
measurement. Linearity and approximation then determine $L_V$
for every admissible source.

We next recover the geometry from $L_V$.
Taking the Laplace transform in time determines the localized
resolvent of $(-\Delta_g)^\alpha$.
Stone's formula then recovers its compressed spectral measure
on $V$, from which we obtain the corresponding compressed
spectral measure of $-\Delta_g$.
This measure determines the local source-to-solution operator
for the wave equation.
The geometric determination result of
\cite[Theorem~2]{Correlation} therefore yields the desired
Riemannian isometry.
\subsection{Motivation and previous literature}
\label{subsec:motivation-literature}

Fractional diffusion equations provide a framework for describing
anomalous transport in media where classical diffusion models do not
adequately capture memory effects or nonlocal spatial interactions.
Their connections with continuous-time random walks and their
applications across physical and biological systems have been
extensively studied \cite{MLP2001,MK2000,MK2004}.
Temporal fractional derivatives account for memory, while spatial
fractional operators describe transport involving long-range
interactions. These models naturally motivate inverse problems:
to what extent can observations of a diffusion process reveal the
properties, and in particular the geometry, of the underlying medium?

A central question in this context concerns the amount of information
needed for unique recovery. For time-fractional diffusion equations,
geometric and coefficient recovery have been established using
restricted boundary observations, including fixed-time measurements
with varying inputs and observations generated by a single suitably
chosen input \cite{KOSY2018,KLLY2021,Kian2022}.
The construction of sources that encode geometric information in one
experiment also appears in inverse problems for the wave equation
\cite{HLO2012}. Most directly related to our work,
Helin--Lassas--Ylinen--Zhang \cite{HLYZ2020} proved that a specially
constructed interior source and the corresponding finite-time local
observation determine a closed Riemannian manifold. Their source
construction provides the starting point for our single-measurement
argument.

A complementary direction is the study of inverse problems for
stationary nonlocal equations. Foundational results for the fractional
Schr\"odinger equation establish potential recovery from partial
exterior data, including uniqueness and reconstruction from a single
measurement \cite{GSU2020,GRSU2020}. On Riemannian manifolds,
Feizmohammadi--Ghosh--Krupchyk--Uhlmann and Choulli--Ouhabaz recover
closed and complete manifolds, respectively, from local
source-to-solution operators associated with the fractional Laplacian,
assuming the geometry is known on the observation region
\cite{FGKU2025,CO2024}. These results demonstrate the geometric
information contained in local fractional measurements and provide
a natural comparison with our finite-time evolution data.

Our approach also connects with inverse spectral theory and the
boundary control method \cite{BK1992,KKL2001}. In particular, the
deterministic theorem of Helin--Lassas--Oksanen--Saksala
\cite{HLOS2018} recovers a complete manifold from its local wave
source-to-solution operator and supplies the final geometric step
of our argument. We extend the single-measurement framework of
\cite{HLYZ2020} to complete, possibly noncompact manifolds, where the
Laplace--Beltrami operator may have continuous spectrum.
A spectral-measure formulation replaces discrete eigenfunction
expansions and enables the reduction to wave data. We thereby address
whether one suitably chosen interior source, observed locally over
a finite time interval, suffices to determine the complete manifold
up to isometry.

\subsection{Organization of the paper}
\label{subsec:organization}

Section~\ref{sec:preliminaries} introduces the notation and
collects the properties of fractional operators and
Mittag--Leffler functions used in the
paper. In Section~\ref{sec:direct-problem}, we establish
well-posedness of the direct problem and define the local
source-to-solution operator.
Section~\ref{sec:inverse-problem} contains the proof of
Proposition \ref{prop:finite-to-infinite-time} and Theorem~\ref{thm:single-measurement}.

\section{Preliminaries}
\label{sec:preliminaries}

We collect the notation and auxiliary facts concerning fractional
operators and Mittag--Leffler functions that
will be used throughout the paper.

\subsection{The fractional Laplace--Beltrami operator}
\label{subsec:prelim-fractional-laplacian}

Let $(M,g)$ be a smooth, connected, complete Riemannian manifold
without boundary, of dimension $n\geq2$. Throughout the paper,
$L^2(M)$ denotes $L^2(M,dV_g)$, where $dV_g$ is the Riemannian
volume measure, and its inner product is taken to be linear
in the first argument.

We use the sign convention that $-\Delta_g$ is nonnegative.
Since $(M,g)$ is complete, the operator $-\Delta_g$, initially
defined on $C_c^\infty(M)$, is essentially self-adjoint on
$L^2(M)$; see \cite{Strichartz1983}. We use the same notation
$-\Delta_g$ for its nonnegative self-adjoint closure.

We recall the spectral formulation used to define its fractional
powers. By the spectral theorem, every nonnegative self-adjoint
operator $A$ on $L^2(M)$ admits a unique projection-valued
measure $P_A$, supported in $[0,\infty)$, such that
\[
    A=\int_{[0,\infty)}\lambda\,dP_A(\lambda).
\]
For each Borel set $B\subset\mathbb R$, the operator $P_A(B)$
is an orthogonal projection on $L^2(M)$. Moreover,
$P_A(\mathbb R)=I$, and $P_A$ is countably additive in the
strong operator topology. We call $P_A$ the spectral measure
of $A$.

Pairing these projections with $h,v\in L^2(M)$ defines the
scalar spectral measure
\[
    \nu_{h,v}^A(B)
    :=
    \langle P_A(B)h,v\rangle_{L^2(M)},
    \qquad B\subset\mathbb R \text{ Borel}.
\]
This is a finite complex measure supported in $[0,\infty)$,
with total variation bounded by
\[
    |\nu_{h,v}^A|(\mathbb R)
    \leq
    \|h\|_{L^2(M)}\|v\|_{L^2(M)}.
\]
In particular, $\nu_{h,h}^A$ is a positive measure, since
\[
    \nu_{h,h}^A(B)
    =
    \|P_A(B)h\|_{L^2(M)}^2.
\]

For a bounded Borel function
$F:[0,\infty)\to\mathbb C$, the spectral functional calculus
defines a bounded operator
\[
    F(A)
    :=
    \int_{[0,\infty)}F(\lambda)\,dP_A(\lambda).
\]
Its matrix elements are given by integration against the
scalar spectral measures:
\begin{equation}\label{eq:prelim-functional-calculus}
    \begin{aligned}
        \langle F(A)h,v\rangle_{L^2(M)}
        &=
        \int_{[0,\infty)}
            F(\lambda)\,d\nu_{h,v}^A(\lambda),\\
        \|F(A)\|_{\mathcal L(L^2(M))}
        &\leq
        \sup_{\lambda\geq0}|F(\lambda)|.
    \end{aligned}
\end{equation}

For $0<\alpha\leq1$, we define
\begin{equation}\label{eq:prelim-fractional-operator}
    (-\Delta_g)^\alpha
    :=
    \int_{[0,\infty)}
        \lambda^\alpha\,dE_g(\lambda),
\end{equation}
with domain
\[
    \mathcal D\bigl((-\Delta_g)^\alpha\bigr)
    :=
    \left\{
        h\in L^2(M):
        \int_{[0,\infty)}
            \lambda^{2\alpha}\,
            d\nu_{h,h}^{-\Delta_g}(\lambda)<\infty
    \right\}.
\]
The domain condition ensures that the spectral integral
defining $(-\Delta_g)^\alpha h$ converges in $L^2(M)$.
The resulting operator is nonnegative and self-adjoint,
and for $\alpha=1$ it coincides with $-\Delta_g$.

The spectral measure of $(-\Delta_g)^\alpha$ is obtained
from $P_{-\Delta_g}$ by the change of spectral variable
$\lambda\mapsto\lambda^\alpha$. More precisely,
\begin{equation}\label{eq:prelim-spectral-pushforward}
    P_{(-\Delta_g)^\alpha}(B)
    =
    P_{-\Delta_g}\bigl(
        \{\lambda\geq0:\lambda^\alpha\in B\}
    \bigr),
    \qquad B\subset\mathbb R \text{ Borel}.
\end{equation}
These definitions apply without any assumption that the
spectrum is discrete.
\subsection{The Caputo derivative and the Laplace transform}
\label{subsec:prelim-caputo}

Throughout the paper, $\beta\in(0,1]$ denotes the temporal
fractional order. Let $X$ be a complex Hilbert space. For
$0<\beta<1$ and
$u\in W_{\mathrm{loc}}^{1,1}([0,\infty);X)$, the Caputo derivative
is defined, for almost every $t>0$, by
\begin{equation}
\label{eq:prelim-caputo}
    \partial_t^\beta u(t)
    :=
    \frac{1}{\Gamma(1-\beta)}
    \int_0^t (t-\tau)^{-\beta}u'(\tau)\,d\tau.
\end{equation}
For $\beta=1$, we set $\partial_t^\beta u=u'$.
All integrals of $X$-valued functions are understood as
Bochner integrals.

For a strongly measurable function $v:[0,\infty)\to X$, its
Laplace transform is
\[
    \widehat v(s)
    :=
    \int_0^\infty e^{-st}v(t)\,dt
\]
whenever the integral converges absolutely in $X$.

Suppose that
$u\in W_{\mathrm{loc}}^{1,1}([0,\infty);X)$ and that, for some
$\sigma\geq0$,
\[
    \int_0^\infty e^{-\sigma t}
    \bigl(\|u(t)\|_X+\|u'(t)\|_X\bigr)\,dt<\infty.
\]
Then, for $\operatorname{Re}s>\sigma$,
\begin{equation}
\label{eq:prelim-caputo-laplace}
    \widehat{\partial_t^\beta u}(s)
    =
    s^\beta\widehat u(s)-s^{\beta-1}u(0),
\end{equation}
where we use the principal branch of the logarithm and set
$ s^\gamma:=\exp(\gamma\operatorname{Log}s)$, $
    s\in\mathbb C_+
    :=
    \{s\in\mathbb C:\operatorname{Re}s>0\}$.
For $0<\beta<1$, this follows by applying Fubini's theorem to
\eqref{eq:prelim-caputo} and using
\[
    \int_0^\infty
    e^{-st}\frac{t^{-\beta}}{\Gamma(1-\beta)}\,dt
    =
    s^{\beta-1}.
\]
For $\beta=1$, it is the usual differentiation formula for the
Laplace transform. Further details on the Caputo derivative
may be found in
\cite{KilbasSrivastavaTrujillo2006,Podlubny1999}.

\subsection{Mittag--Leffler functions}
\label{subsec:prelim-mittag-leffler}

For $\rho,\eta>0$, the two-parameter Mittag--Leffler function is
defined by
\begin{equation}
\label{eq:prelim-ML-definition}
    E_{\rho,\eta}(z)
    :=
    \sum_{k=0}^\infty
    \frac{z^k}{\Gamma(\rho k+\eta)},
    \qquad z\in\mathbb C.
\end{equation}
The series converges locally uniformly in $\mathbb C$, and
$E_{\rho,\eta}$ is entire. In particular,
\[
    E_{\rho,\eta}(0)=\frac{1}{\Gamma(\eta)},
    \qquad
    E_{1,1}(z)=e^z.
\]
We shall use the functions $E_{\beta,1}$ and
$E_{\beta,\beta}$.

Termwise differentiation of
\eqref{eq:prelim-ML-definition} gives
\begin{equation}
\label{eq:prelim-ML-derivative}
    \frac{d}{dz}E_{\beta,1}(-\lambda z^\beta)
    =
    -\lambda z^{\beta-1}
    E_{\beta,\beta}(-\lambda z^\beta),
    \qquad
    z\in\mathbb C_+,\quad \lambda\geq0,
\end{equation}
where the complex powers are defined using the principal branch, and 
$
    \mathbb C_+:=\{z\in\mathbb C:\operatorname{Re}z>0\}.
$

We record the sectorial estimate used below. Let $0<\beta<1$,
$\eta>0$, and
\[
    \frac{\beta\pi}{2}<\theta<\pi.
\]
Then there exists $C=C(\beta,\eta,\theta)>0$ such that
\begin{equation}\label{eq:prelim-ML-sectorial}
    |E_{\beta,\eta}(\zeta)|
    \leq \frac{C}{1+|\zeta|},
    \qquad
    \zeta\neq0,\quad
    \theta\leq|\arg\zeta|\leq\pi.
\end{equation}
The estimate extends to $\zeta=0$ by continuity.
It follows from the sectorial asymptotic expansion
of the Mittag--Leffler function together with its boundedness
on compact sets; see
\cite{KilbasSrivastavaTrujillo2006,HauboldMathaiSaxena2011}.

Taking $\zeta=-r$ in \eqref{eq:prelim-ML-sectorial} yields
\begin{equation}
\label{eq:prelim-ML-negative-axis}
    |E_{\beta,1}(-r)|
    +
    |E_{\beta,\beta}(-r)|
    \leq \frac{C_\beta}{1+r},
    \qquad r\geq0.
\end{equation}
The same estimate holds for $\beta=1$, since
$E_{1,1}(-r)=e^{-r}$. Consequently, for
$0<\beta\leq1$, $t>0$, and $\lambda\geq0$,
\begin{align}
\label{eq:prelim-ML-real-kernels}
    |E_{\beta,1}(-\lambda t^\beta)|
    &\leq \frac{C_\beta}{1+\lambda t^\beta},
    \\
    |t^{\beta-1}E_{\beta,\beta}(-\lambda t^\beta)|
    &\leq
    \frac{C_\beta t^{\beta-1}}{1+\lambda t^\beta}.
    \notag
\end{align}

We will use the following Laplace-transform identities later.
\begin{align}
\label{eq:prelim-ML-laplace-homogeneous}
    \int_0^\infty
    e^{-st}E_{\beta,1}(-\lambda t^\beta)\,dt
    &=
    \frac{s^{\beta-1}}{s^\beta+\lambda},
    \\
\label{eq:prelim-ML-laplace-inhomogeneous}
    \int_0^\infty
    e^{-st}t^{\beta-1}
    E_{\beta,\beta}(-\lambda t^\beta)\,dt
    &=
    \frac{1}{s^\beta+\lambda},
\end{align}
for $\lambda\geq0$ and $\operatorname{Re}s>0$.
For each fixed $\lambda\geq0$, these identities follow for real
$s>\lambda^{1/\beta}$ by integrating the defining series term
by term. The estimates in
\eqref{eq:prelim-ML-real-kernels} imply that the left-hand sides
are holomorphic functions of $s$ in
$\mathbb C_+:=\{s\in\mathbb C:\operatorname{Re}s>0\}$.
The right-hand sides are also holomorphic there, with
$s^\beta$ defined using the principal branch, since
$s^\beta+\lambda\neq0$.
The identity theorem therefore extends both identities to
every $s\in\mathbb C_+$.

\section{Analysis of the direct problem}\label{sec:direct-problem}

Let $(M,g)$ be a connected, smooth, complete Riemannian manifold
without boundary, and let
\[
    H:=L^2(M,dV_g).
\]
We denote by $L_g$ the nonnegative self-adjoint realization of
$-\Delta_g$ on $H$ and set
\[
    A_g:=L_g^\alpha,
    \qquad 0<\alpha<1.
\]
Throughout, $\alpha$ denotes the spatial fractional order, whereas
$\beta\in(0,1]$ denotes the temporal fractional order.

\subsection{Existence and uniqueness of strong solutions}
\label{subsec:well-posedness}

Fix $\alpha,\beta\in(0,1]$. We consider the initial value problem
\begin{equation}\label{eq:pde}
    \partial_t^\beta u(t)+(-\Delta_g)^\alpha u(t)=f(t),
    \qquad t>0,
\end{equation}
with
\begin{equation}\label{eq:initial}
    u(0)=0.
\end{equation}
Here $\partial_t^\beta$ is the Caputo derivative when
$0<\beta<1$ and the ordinary derivative when $\beta=1$.
We establish existence and uniqueness for the regular sources
used in the inverse problem and derive a representation of
the solution.

Let $P_{(-\Delta_g)^\alpha}$ denote the spectral measure of
$(-\Delta_g)^\alpha$. For $t>0$, define
\[
    S_0(t):=
    E_{\beta,1}\bigl(-t^\beta(-\Delta_g)^\alpha\bigr),
    \qquad
    S_1(t):=
    t^{\beta-1}
    E_{\beta,\beta}\bigl(-t^\beta(-\Delta_g)^\alpha\bigr),
\]
and set $S_0(0)=I$. These operators are defined by the
spectral calculus:
\begin{align*}
    S_0(t)h
    &=
    \int_{[0,\infty)}
        E_{\beta,1}(-\lambda t^\beta)\,
        dP_{(-\Delta_g)^\alpha}(\lambda)h,\\
    S_1(t)h
    &=
    \int_{[0,\infty)}
        t^{\beta-1}E_{\beta,\beta}(-\lambda t^\beta)\,
        dP_{(-\Delta_g)^\alpha}(\lambda)h,
    \qquad h\in L^2(M).
\end{align*}
The real-axis estimates \eqref{eq:prelim-ML-real-kernels} imply
\begin{align}
    \|S_1(t)\|_{\mathcal L(L^2(M))}
    &\leq
    t^{\beta-1}
    \sup_{\lambda\geq0}
        |E_{\beta,\beta}(-\lambda t^\beta)|
    \notag\\
    &\leq C_\beta t^{\beta-1},
    \qquad t>0.
    \label{eq:S1-operator-bound}
\end{align}
Moreover, since
\[
    0\leq E_{\beta,1}(-r)\leq1,
    \qquad r\geq0,
\]
we have
\begin{equation}\label{eq:solution-operator-bounds}
    \|S_0(t)\|_{\mathcal L(L^2(M))}\leq1,
    \qquad
    \|I-S_0(t)\|_{\mathcal L(L^2(M))}\leq1,
    \qquad t\geq0.
\end{equation}
Dominated convergence in the spectral representation shows
that $S_0$ is strongly continuous on $[0,\infty)$ and $S_1$
is strongly continuous on $(0,\infty)$.
For $\beta=1$, both operators coincide with
$e^{-t(-\Delta_g)^\alpha}$.

\begin{defin}[Strong solution]\label{def:strong-solution}
Let $f\in C([0,\infty);L^2(M))$.
A function
\[
    u\in C^1([0,\infty);L^2(M))
\]
is called a strong solution of
\eqref{eq:pde}--\eqref{eq:initial} if $u(0)=0$,
\[
    u(t)\in\mathcal D((-\Delta_g)^\alpha)
    \qquad\text{for every }t\geq0,
\]
and
\[
    \partial_t^\beta u(t)+(-\Delta_g)^\alpha u(t)=f(t)
\]
holds in $L^2(M)$ for every $t>0$.
\end{defin}

In the proof of the following proposition, the operator $S_1$ gives the solution representation, while
$S_0$ will be used in an integration-by-parts argument to
establish spatial regularity.

\begin{prop}\label{prop:strong-existence}
Suppose that
\[
    f\in C_c^2\bigl((0,\infty);L^2(M)\bigr),
\]
with the convention $f(0)=0$.
Then
\begin{equation}\label{eq:strong-solution-representation}
    u(t)
    :=
    \int_0^t
        (t-s)^{\beta-1}
        E_{\beta,\beta}
        \bigl(-(t-s)^\beta(-\Delta_g)^\alpha\bigr)
        f(s)\,ds
\end{equation}
is the unique strong solution of
\eqref{eq:pde}--\eqref{eq:initial}.
Moreover,
\[
    u\in C^1([0,\infty);L^2(M))
    \cap
    C\bigl([0,\infty);
        \mathcal D((-\Delta_g)^\alpha)\bigr),
\]
where $\mathcal D((-\Delta_g)^\alpha)$ is equipped with
its graph norm.
\end{prop}

\begin{proof}
Choose $T>0$ such that $\operatorname{supp}f\subset(0,T)$.
We first verify that the integral defining $u(t)$ exists.
For $t>0$,
\begin{align*}
    \int_0^t\|S_1(t-s)f(s)\|_{L^2(M)}\,ds
    &\leq
    C_\beta
    \sup_{0\leq s\leq t}\|f(s)\|_{L^2(M)}
    \int_0^t(t-s)^{\beta-1}\,ds\\
    &=
    \frac{C_\beta}{\beta}t^\beta
    \sup_{0\leq s\leq t}\|f(s)\|_{L^2(M)}.
\end{align*}
Thus the integrand is Bochner integrable.

To prove time regularity, we differentiate the convolution.
Changing variables gives
\[
    u(t)=\int_0^t S_1(r)f(t-r)\,dr.
\]
The integrable bound on $S_1$ and the smoothness of $f$
justify differentiation under the integral.
Using $f(0)=0$, we obtain
\begin{align}
    u'(t)
    &=
    S_1(t)f(0)
    +\int_0^tS_1(r)f'(t-r)\,dr
    \notag\\
    &=
    \int_0^tS_1(t-s)f'(s)\,ds,
    \qquad t>0.
    \label{eq:strong-solution-time-derivative}
\end{align}
Dominated convergence gives continuity of $u'$ for $t>0$.
Since $f$ vanishes near $t=0$, so does $u$. Hence
\[
    u\in C^1([0,\infty);L^2(M)),
    \qquad u(0)=0.
\]

To establish spatial regularity, we first restrict the
calculation to bounded spectral subspaces.
For $N\geq1$, let
\[
    P_N:=P_{(-\Delta_g)^\alpha}([0,N]).
\]
These projections satisfy
$P_NL^2(M)\subset D((-\Delta_g)^\alpha)$,
$\|(-\Delta_g)^\alpha P_N\|\leq N$, and $P_N\to I$ strongly.
Set
\[
    u_N(t):=P_Nu(t)
    =
    \int_0^t S_1(t-s)P_Nf(s)\,ds.
\]
Boundedness of $P_N$ gives
\[
    u_N\in C^1([0,\infty);L^2(M)),
    \qquad u_N'=P_Nu'.
\]

Integration by parts will give an expression for
$(-\Delta_g)^\alpha u_N$ in terms of $f'$.
The Mittag--Leffler differentiation identity reads
\[
    \frac{d}{ds}\bigl(S_0(t-s)P_N\bigr)
    =
    (-\Delta_g)^\alpha S_1(t-s)P_N,
    \qquad 0<s<t.
\]
The bound
\[
    \|(-\Delta_g)^\alpha S_1(t-s)P_N\|_{\mathcal L(L^2(M))}
    \leq NC_\beta(t-s)^{\beta-1}
\]
shows that this derivative is integrable up to $s=t$.
Integrating by parts and using $f(0)=0$, we obtain
\begin{align}
    (-\Delta_g)^\alpha u_N(t)
    &=
    P_Nf(t)
    -
    \int_0^tS_0(t-s)P_Nf'(s)\,ds
    \notag\\
    &=
    \int_0^t
    \bigl(I-S_0(t-s)\bigr)P_Nf'(s)\,ds.
    \label{eq:truncated-spatial-identity}
\end{align}

We use this identity to prove uniform convergence of
$(-\Delta_g)^\alpha u_N$ on $[0,\infty)$.
For $K>N$, the Cauchy--Schwarz inequality and the bound
$\|I-S_0(t)\|\leq1$ give
\begin{align}
    &\sup_{t\geq0}
    \|(-\Delta_g)^\alpha u_K(t)
      -(-\Delta_g)^\alpha u_N(t)\|_{L^2(M)}^2
    \notag\\
    &\qquad\leq
    T\int_0^T
    \|(P_K-P_N)f'(s)\|_{L^2(M)}^2\,ds
    \notag\\
    &\qquad\leq
    T\int_0^T
    \|(I-P_N)f'(s)\|_{L^2(M)}^2\,ds.
    \label{eq:spectral-tail-A}
\end{align}
Here we used that $f'$ vanishes outside $(0,T)$.
Strong convergence of $P_N$ and dominated convergence
show that the last integral tends to zero.
Therefore $(-\Delta_g)^\alpha u_N$ converges uniformly
on $[0,\infty)$ to an $L^2(M)$-valued function $w$.

We next verify continuity of this limit.
Each function
$(-\Delta_g)^\alpha u_N=(-\Delta_g)^\alpha P_Nu$
is continuous because $(-\Delta_g)^\alpha P_N$ is bounded
and $u$ is continuous.
Uniform convergence then gives
\[
    w\in C([0,\infty);L^2(M)).
\]

The closedness of $(-\Delta_g)^\alpha$ allows us to identify
$w$ with $(-\Delta_g)^\alpha u$.
For every fixed $t\geq0$,
\[
    u_N(t)=P_Nu(t)\longrightarrow u(t),
    \qquad
    (-\Delta_g)^\alpha u_N(t)\longrightarrow w(t)
    \quad\text{in }L^2(M).
\]
Since $(-\Delta_g)^\alpha$ is self-adjoint, it is closed,
and these limits imply
\[
    u(t)\in D((-\Delta_g)^\alpha),
    \qquad
    (-\Delta_g)^\alpha u(t)=w(t).
\]

We finish the spatial regularity argument by checking
continuity in the graph norm.
For every $t_0\geq0$,
\begin{align*}
    \|u(t)-u(t_0)\|_{D((-\Delta_g)^\alpha)}^2
    &=
    \|u(t)-u(t_0)\|_{L^2(M)}^2
    +
    \|w(t)-w(t_0)\|_{L^2(M)}^2\\
    &\longrightarrow0
    \qquad\text{as }t\to t_0,
\end{align*}
because $u$ and $w$ are continuous in $L^2(M)$.
Thus
\[
    u\in C([0,\infty);D((-\Delta_g)^\alpha)).
\]

We now verify the fractional diffusion equation on each
bounded spectral subspace.
On $P_NL^2(M)$, the operator $(-\Delta_g)^\alpha$ is bounded
and the Mittag--Leffler series converge in operator norm.
Together with the integrable bound on the time kernel,
this permits us to apply the scalar solution formula in
\cite[Proposition~5]{HelinLassasYlinenZhang2019}
through spectral calculus on $[0,N]$.
We obtain
\begin{equation}\label{eq:truncated-fractional-equation}
    \partial_t^\beta u_N(t)
    +(-\Delta_g)^\alpha u_N(t)=P_Nf(t).
\end{equation}

To recover the equation for $u$, we commute $P_N$ with
the Caputo derivative and pass to the limit.
For $0<\beta<1$, boundedness of $P_N$ gives
\begin{align*}
    \partial_t^\beta u_N(t)
    &=
    \frac1{\Gamma(1-\beta)}
    \int_0^t(t-s)^{-\beta}P_Nu'(s)\,ds\\
    &=
    P_N\partial_t^\beta u(t).
\end{align*}
For $\beta=1$, the same identity follows from
$u_N'=P_Nu'$.
Letting $N\to\infty$ in
\eqref{eq:truncated-fractional-equation} yields
\[
    \partial_t^\beta u(t)+(-\Delta_g)^\alpha u(t)=f(t).
\]

For uniqueness, we reduce the difference of two strong
solutions to a homogeneous equation on each spectral
subspace.
Let $v$ denote this difference and set $v_N=P_Nv$.
Projecting the equation gives
\[
    \partial_t^\beta v_N(t)+(-\Delta_g)^\alpha v_N(t)=0,
    \qquad v_N(0)=0.
\]

We prove that $v_N$ vanishes by iterating a fractional
integral estimate.
Applying the fractional integral of order $\beta$, with
ordinary integration when $\beta=1$, gives
\[
    v_N(t)
    =
    -\frac1{\Gamma(\beta)}
    \int_0^t(t-s)^{\beta-1}
    (-\Delta_g)^\alpha v_N(s)\,ds.
\]
Since $v_N(s)\in P_NL^2(M)$,
\[
    \|v_N(t)\|_{L^2(M)}
    \leq
    \frac{N}{\Gamma(\beta)}
    \int_0^t(t-s)^{\beta-1}\|v_N(s)\|_{L^2(M)}\,ds.
\]
Iterating this inequality gives, for every integer $m\geq1$,
\[
    \|v_N(t)\|_{L^2(M)}
    \leq
    \frac{N^m t^{m\beta}}{\Gamma(m\beta+1)}
    \sup_{0\leq s\leq t}\|v_N(s)\|_{L^2(M)}.
\]
For fixed $N$ and $t$, the supremum is finite by continuity
of $v_N$, while the coefficient tends to zero as
$m\to\infty$.
Thus $v_N(t)=0$. Passing to the strong limit $P_N\to I$
gives $v(t)=0$, proving uniqueness.
\end{proof}

\begin{rmk}\label{rem:less-regular-sources}
The representation
\eqref{eq:strong-solution-representation} also extends to
less regular sources. Let
$f\in L^p((0,T);L^2(M))$, extended by zero outside $(0,T)$,
where $p=1$ if $\beta=1$ and $1/\beta<p\leq\infty$
if $0<\beta<1$.
The same formula defines a function
$u\in C([0,\infty);L^2(M))$ with $u(0)=0$, usually called
the mild solution. For $0<\beta<1$, it satisfies
\begin{equation}\label{eq:mild-solution-estimate}
    \|u(t)\|_{L^2(M)}
    \leq
    C_{\beta,p}t^{\beta-1/p}
    \|f\|_{L^p((0,t);L^2(M))},
    \qquad t>0,
\end{equation}
with $1/p=0$ when $p=\infty$.
For $\beta=1$, the corresponding bound is
$\|u(t)\|_{L^2(M)}
\leq\|f\|_{L^1((0,t);L^2(M))}$.
These estimates follow from the bounds on $S_1$ and
H\"older's inequality. Continuity follows by approximation
with smooth sources; when $p=\infty$, the approximation
is taken in $L^q$ on bounded time intervals for a finite
$q>1/\beta$.

For $t>T$, the representation reduces to
\begin{equation}\label{eq:mild-after-T}
    u(t)=\int_0^T S_1(t-s)f(s)\,ds.
\end{equation}

For $t>T$, this integral is well defined even under the weaker
assumption $f\in L^1((0,T);L^2(M))$, because
$t-s\in[t-T,t]$ stays away from zero.
The argument establishing Bochner integrability therefore
remains valid, and strong continuity of $S_1$, together with
dominated convergence, gives continuity of $u$ on $(T,\infty)$.
Likewise, the time-analyticity argument below uses holomorphy
and uniform bounds of the kernel on compact subsets of the
complex-time domain and does not require derivatives of $f$.

For convenience, we still work with strong solutions throughout the inverse
problem, so that the equation holds in $L^2(M)$ at every
 positive time.

\end{rmk}

\subsection{The local source-to-solution operator $L_V$}
\label{subsec:source-to-solution}

We now verify that the local source-to-solution operator
introduced in the introduction is well defined and record
its relation to the localized resolvent.
Recall that $e_V$ denotes extension by zero from $V$ to $M$
and $r_V=e_V^*$ denotes restriction to $V$.

For $f\in C_c^2((0,\infty);L^2(V))$, boundedness of $e_V$
implies
\[
    e_Vf\in C_c^2((0,\infty);L^2(M)).
\]
Proposition~\ref{prop:strong-existence} therefore gives a
unique strong solution $u_f$ of
\begin{equation}\label{eq:localized-source-problem}
    \partial_t^\beta u_f+(-\Delta_g)^\alpha u_f=e_Vf,
    \qquad t>0,
    \qquad u_f(0)=0.
\end{equation}
Consequently, the operator $L_Vf=r_Vu_f$ has the mapping
property
\[
    L_V:
    C_c^2((0,\infty);L^2(V))
    \longrightarrow C^1([0,\infty);L^2(V)).
\]
The representation
\eqref{eq:strong-solution-representation} gives
\begin{equation}\label{eq:local-source-to-solution}
    (L_Vf)(t)
    =
    \int_0^t r_VS_1(t-s)e_Vf(s)\,ds.
\end{equation}
In particular, $L_V$ is linear and causal: its value at time
$t$ depends only on the source on $[0,t]$.
Since $\|e_V\|=\|r_V\|=1$, the bound
\eqref{eq:S1-operator-bound} and the derivative formula
\eqref{eq:strong-solution-time-derivative} imply that, for
every $R>0$ and $j=0,1$,
\begin{equation}\label{eq:local-source-to-solution-bound}
    \sup_{0\leq t\leq R}
    \|\partial_t^j(L_Vf)(t)\|_{L^2(V)}
    \leq
    \frac{C_\beta}{\beta}R^\beta
    \sup_{0\leq t\leq R}
    \|f^{(j)}(t)\|_{L^2(V)}.
\end{equation}

For $T'>0$, we write
\[
    L_{V,T'}f:=(L_Vf)|_{(0,T')},
    \qquad
    f\in C_c^2((0,T');L^2(V)),
\]
for the finite-time source-to-solution operator, where
$f$ is extended by zero to $(0,\infty)$.

The associated localized resolvent is
\begin{equation}\label{eq:localized-resolvent}
    L_V(z)
    :=
    r_V\bigl((-\Delta_g)^\alpha-z\bigr)^{-1}e_V,
    \qquad
    z\in\mathbb C\setminus[0,\infty).
\end{equation}
Its relation to the time-dependent operator is obtained
by taking the Laplace transform of
\eqref{eq:local-source-to-solution}.
For
\[
    \widehat f(s):=\int_0^\infty e^{-st}f(t)\,dt,
    \qquad s>0,
\]
we have
\begin{equation}\label{eq:laplace-source-to-solution}
    \int_0^\infty e^{-st}(L_Vf)(t)\,dt
    =
    L_V(-s^\beta)\widehat f(s),
    \qquad s>0.
\end{equation}
Indeed, the scalar identity
\[
    \int_0^\infty
        e^{-st}t^{\beta-1}
        E_{\beta,\beta}(-\lambda t^\beta)\,dt
    =
    \frac{1}{s^\beta+\lambda},
    \qquad s>0,\quad\lambda\geq0,
\]
and the spectral calculus give
\[
    \int_0^\infty e^{-st}S_1(t)\,dt
    =
    \bigl((-\Delta_g)^\alpha+s^\beta\bigr)^{-1},
\]
where the operator integral is understood in the strong
sense. Moreover,
\begin{align*}
    &\int_0^\infty\int_0^t
        e^{-st}\|S_1(t-\tau)e_Vf(\tau)\|_{L^2(M)}
        \,d\tau\,dt\\
    &\qquad\leq
    C_\beta\Gamma(\beta)s^{-\beta}
    \int_0^\infty
        e^{-s\tau}\|f(\tau)\|_{L^2(V)}\,d\tau
    <\infty.
\end{align*}
Thus Fubini's theorem applies to the convolution in
\eqref{eq:local-source-to-solution} and yields
\eqref{eq:laplace-source-to-solution}.

\section{Analysis of the inverse problem}\label{sec:inverse-problem}

\subsection{Time analyticity after the source vanishes}
\label{subsec:time-analyticity}

We now show that the solution is analytic in time beyond the
temporal support of the source. The argument uses only the
nonnegativity and self-adjointness of $A_g$.

Let $\mathbb C_+:=\{z\in\mathbb C:\operatorname{Re}z>0\}$,
and define complex powers on $\mathbb C_+$ using the principal
branch of the logarithm.

The following lemmas establish uniform kernel bounds and
holomorphy of the solution operators in complex time.
\begin{lemma}\label{lem:uniform-bounds}
For every compact set $K\Subset\mathbb C_+$, there exists $C_K>0$
such that
\begin{align}
    |E_{\beta,1}(-\lambda z^\beta)|
    &\leq C_K,\label{eq:S0-bound}\\
    |z^{\beta-1}E_{\beta,\beta}(-\lambda z^\beta)|
    &\leq C_K,\label{eq:S1-bound}
\end{align}
for all $z\in K$ and $\lambda\geq0$.
\end{lemma}

\begin{proof}
For $\beta=1$, both estimates follow from
\[
    |e^{-\lambda z}|
    =
    e^{-\lambda\operatorname{Re}z}
    \leq1,
    \qquad z\in\mathbb C_+,\quad \lambda\geq0.
\]

Suppose that $0<\beta<1$. Choose $\theta$ such that
\[
    \frac{\beta\pi}{2}
    <
    \theta
    <
    \pi-\frac{\beta\pi}{2}.
\]
Such a choice is possible because $\beta<1$.
For $z\in\mathbb C_+$ and $\lambda>0$, the principal branch
of the complex power gives
\[
    |\arg(-\lambda z^\beta)|
    =
    \pi-\beta|\arg z|
    >
    \pi-\frac{\beta\pi}{2}
    >
    \theta.
\]
Applying \eqref{eq:prelim-ML-sectorial} with
$\zeta=-\lambda z^\beta$ and $\eta\in\{1,\beta\}$ yields
\[
    |E_{\beta,\eta}(-\lambda z^\beta)|
    \leq
    \frac{C}{1+\lambda|z|^\beta}
    \leq C.
\]
The same bound holds for $\lambda=0$, after increasing $C$
if necessary, since
\[
    E_{\beta,\eta}(0)=\frac{1}{\Gamma(\eta)}.
\]

Since $K$ is compact and contained in $\mathbb C_+$,
\[
    r_K:=\min_{z\in K}|z|>0.
\]
Consequently,
\[
    |z^{\beta-1}|
    =
    |z|^{\beta-1}
    \leq r_K^{\beta-1},
    \qquad z\in K.
\]
The two asserted bounds follow with
\[
    C_K=C\max\{1,r_K^{\beta-1}\}.
\]
\end{proof}

For $z\in\mathbb C_+$, we may therefore define
\[
    S_0(z):=
    E_{\beta,1}\bigl(-z^\beta(-\Delta_g)^\alpha\bigr),
    \qquad
    S_1(z):=
    z^{\beta-1}
    E_{\beta,\beta}\bigl(-z^\beta(-\Delta_g)^\alpha\bigr)
\]
as bounded operators on $L^2(M)$.
Both families are locally uniformly bounded in
$\mathcal L(L^2(M))$.

\begin{lemma}\label{lem:holomorphic-multipliers}
For each $h\in L^2(M)$, the functions
\[
    z\longmapsto S_0(z)h,
    \qquad
    z\longmapsto S_1(z)h
\]
are holomorphic on $\mathbb C_+$ with values in $L^2(M)$.
\end{lemma}

\begin{proof}
We prove the assertion for $S_1$; the proof for $S_0$ is analogous.
For $h,v\in L^2(M)$, let
\[
    \nu_{h,v}^{(-\Delta_g)^\alpha}(B)
    :=
    \langle P_{(-\Delta_g)^\alpha}(B)h,v\rangle_{L^2(M)}
\]
be the scalar spectral measure for $(-\Delta_g)^\alpha$.
This is a finite complex measure satisfying
\[
    \bigl|\nu_{h,v}^{(-\Delta_g)^\alpha}\bigr|([0,\infty))
    \leq
    \|h\|_{L^2(M)}\|v\|_{L^2(M)}.
\]
The spectral calculus gives
\[
    \langle S_1(z)h,v\rangle_{L^2(M)}
    =
    \int_{[0,\infty)}
        z^{\beta-1}E_{\beta,\beta}(-\lambda z^\beta)\,
        d\nu_{h,v}^{(-\Delta_g)^\alpha}(\lambda).
\]
For every $\lambda\geq0$, the integrand is holomorphic in
$z\in\mathbb C_+$.
Lemma~\ref{lem:uniform-bounds} provides a bound uniform in
$\lambda$ on compact subsets of $\mathbb C_+$.
Set
\[
    F(z):=\langle S_1(z)h,v\rangle_{L^2(M)}.
\]
Since $\nu_{h,v}^{(-\Delta_g)^\alpha}$ has finite total variation,
the preceding bound on compact subsets and dominated convergence
show that $F$ is continuous on $\mathbb C_+$.

Let $\mathcal T$ be a triangle with
$\overline{\mathcal T}\subset\mathbb C_+$.
Lemma~\ref{lem:uniform-bounds} gives
\begin{align*}
    &\int_{\partial\mathcal T}\int_{[0,\infty)}
    \left|z^{\beta-1}E_{\beta,\beta}(-\lambda z^\beta)\right|
    \,d\bigl|\nu_{h,v}^{(-\Delta_g)^\alpha}\bigr|(\lambda)
    \,|dz|
    \leq
    C_{\mathcal T}
    \bigl|\nu_{h,v}^{(-\Delta_g)^\alpha}\bigr|([0,\infty))
    <\infty.
\end{align*}
Fubini's theorem therefore permits us to interchange the
contour integral and the spectral integral:
\[
    \int_{\partial\mathcal T}F(z)\,dz
    =
    \int_{[0,\infty)}
    \left(
        \int_{\partial\mathcal T}
        z^{\beta-1}E_{\beta,\beta}(-\lambda z^\beta)\,dz
    \right)
    d\nu_{h,v}^{(-\Delta_g)^\alpha}(\lambda)
    =0.
\]
The inner integral vanishes by Cauchy's theorem, since its
integrand is holomorphic in $z\in\mathbb C_+$ for each fixed
$\lambda\geq0$.
Thus $F$ is continuous and its integral around every such
triangle is zero. Morera's theorem implies that $F$ is
holomorphic on $\mathbb C_+$.


Since $v\in L^2(M)$ was arbitrary, we have $z\mapsto S_1(z)h$ is weakly holomorphic.
The operator bounds established above also show that this
map is locally bounded in $L^2(M)$.
The standard equivalence between weak and strong holomorphy
for locally bounded Banach-space-valued maps therefore implies
that $z\mapsto S_1(z)h$ is holomorphic in the $L^2(M)$ norm.
In particular, its complex derivative exists as a limit in
$L^2(M)$, as required.
\end{proof}

We have the following results for the holomorphic extension of $L_V f$.

\begin{thm}\label{thm:main}
Let $u$ be the strong solution of
\eqref{eq:pde}--\eqref{eq:initial} with
$f\in C_c^2((0,\infty);L^2(M))$, and suppose that $f(t)=0$
for $t>T$, where $T>0$.
Then $u$ is real analytic on $(T,\infty)$ with values in
$L^2(M)$.
More precisely, it admits a holomorphic extension to
\begin{equation}\label{eq:OmegaT}
    \Omega_T:=\{z\in\mathbb C:\operatorname{Re}z>T\},
\end{equation}
given by
\begin{equation}\label{eq:holomorphic-extension}
    \begin{aligned}
        U(z)
        &:=
        \int_0^T S_1(z-s)f(s)\,ds\\
        &=
        \int_0^T
        (z-s)^{\beta-1}
        E_{\beta,\beta}
        \bigl(-(z-s)^\beta(-\Delta_g)^\alpha\bigr)f(s)\,ds.
    \end{aligned}
\end{equation}
In particular, if
$f\in C_c^2((0,\infty);L^2(V))$ vanishes for $t>T$,
then $L_Vf$ is real analytic on $(T,\infty)$ with values in
$L^2(V)$.
\end{thm}

\begin{proof}
Fix a compact set $K\Subset\Omega_T$. Then
\[
    K_T:=\{z-s:z\in K,\ 0\leq s\leq T\}
\]
is a compact subset of $\mathbb C_+$.
Lemma~\ref{lem:uniform-bounds} gives
\[
    \|S_1(z-s)f(s)\|_{L^2(M)}
    \leq C_K\|f(s)\|_{L^2(M)},
    \qquad z\in K,\quad s\in(0,T).
\]
Since $f\in L^1((0,T);L^2(M))$, the integral in
\eqref{eq:holomorphic-extension} is well defined and locally
uniformly bounded.


The preceding estimate and dominated convergence show that
$U$ is continuous with values in $L^2(M)$.
Fix $v\in L^2(M)$ and consider
\[
    F_v(z):=\langle U(z),v\rangle_{L^2(M)}
    =
    \int_0^T
    \langle S_1(z-s)f(s),v\rangle_{L^2(M)}\,ds.
\]
For each fixed $s\in(0,T)$, the integrand is holomorphic in
$z$ by Lemma~\ref{lem:holomorphic-multipliers}.
If $\mathcal T$ is a triangle whose closure is contained in
$\Omega_T$, the same estimate permits us to apply Fubini's
theorem and obtain
\[
    \int_{\partial\mathcal T}F_v(z)\,dz
    =
    \int_0^T
    \left(
        \int_{\partial\mathcal T}
        \langle S_1(z-s)f(s),v\rangle_{L^2(M)}\,dz
    \right)ds
    =0.
\]
Morera's theorem therefore shows that $F_v$ is holomorphic, which says $U$ is weakly holomorphic. Local boundedness then implies
$L^2(M)$-valued holomorphy.



By \eqref{eq:strong-solution-representation} and the assumption
that $f$ vanishes for $t>T$, we have $U(t)=u(t)$ for $t>T$.

For a source $f\in C_c^2((0,\infty);L^2(V))$ vanishing
for $t>T$, apply the preceding argument to the strong solution
with source $e_Vf$.
Since $r_V:L^2(M)\to L^2(V)$ is bounded, restriction of its
holomorphic extension to $V$ gives a holomorphic extension
of $L_Vf$ to $\Omega_T$.
\end{proof}

\begin{rmk}[The Euclidean representation]

When $M=\mathbb R^n$ with its Euclidean metric,
$(-\Delta)^\alpha$ has Fourier symbol $|\xi|^{2\alpha}$.
For $t>T$, the spatial Fourier transform of the solution is
\begin{equation}\label{eq:Fourier-real}
\begin{aligned}
    \widehat u(\xi,t)
    &=
    E_{\beta,1}(-t^\beta|\xi|^{2\alpha})\widehat u_0(\xi)\\
    &\quad+
    \int_0^T
    (t-s)^{\beta-1}
    E_{\beta,\beta}(-(t-s)^\beta|\xi|^{2\alpha})
    \widehat f(\xi,s)\,ds.
\end{aligned}
\end{equation}
The holomorphic extension satisfies
\begin{equation}\label{eq:Fourier-complex}
\begin{aligned}
    \widehat U(\xi,z)
    &=
    E_{\beta,1}(-z^\beta|\xi|^{2\alpha})\widehat u_0(\xi)\\
    &\quad+
    \int_0^T
    (z-s)^{\beta-1}
    E_{\beta,\beta}(-(z-s)^\beta|\xi|^{2\alpha})
    \widehat f(\xi,s)\,ds,
    \qquad z\in\Omega_T.
\end{aligned}
\end{equation}
These integrals are understood in $L^2(\mathbb R^n_\xi)$.
Plancherel's theorem and Lemma~\ref{lem:uniform-bounds} give the
$L^2$-valued holomorphy directly from this representation.
\end{rmk}

Now, we are ready to use the analyticity established to continue finite-time
observations beyond the support of the source.
\begin{prop}
Let $T'>0$ and
$f\in C_c^2((0,T');L^2(V))$.
If
\[
    L_{V,T'}f=\widetilde L_{\widetilde V,T'}f,
\]
then
\[
    L_Vf=\widetilde L_{\widetilde V}f.
\]
\end{prop}

\begin{proof}
The compactness of the temporal support of $f$ gives a number
$T_0<T'$ such that $f(t)=0$ for $t>T_0$.
By Theorem~\ref{thm:main}, both $L_Vf$ and
$\widetilde L_{\widetilde V}f$ are real analytic on
$(T_0,\infty)$ with values in $L^2(V)$.
Their difference vanishes on $(T_0,T')$.
Applying the identity theorem after pairing with an arbitrary
element of $L^2(V)$ shows that it vanishes on $(T_0,\infty)$.
Together with the assumed equality on $(0,T')$, this proves
the assertion.
\end{proof}

The time analyticity established above is the main ingredient
in recovering the full local source-to-solution operator from
a single measurement. Following the source construction in
\cite{HelinLassasYlinenZhang2019}, we can choose $h$ as a countable
superposition of separated temporal pulses with spatial
profiles drawn from a dense family in $L^2(V)$.
Causality ensures that later pulses do not affect the
observations made before they begin. After subtracting the
responses to the preceding pulses, we can therefore isolate
the response to each pulse during its support and in the gap
that follows. Since this response is analytic after the pulse
ends, its values in that gap determine its values at all later
times. Repeating this argument recovers the complete response
to every pulse from the single finite-time measurement.

The temporal profiles are chosen so that translated and
normalized pulses can be used to approximate general time
profiles. Together with the density of the spatial profiles,
linearity, and the continuity estimates for $L_V$, this
determines the response to every admissible source.

\subsection{Determination of the manifold from $L_V$}
\label{subsec:manifold-determination}

We now recover the local spectral measure from the
source-to-solution map and reduce the geometric inverse problem
to the corresponding problem for the wave equation to determine
the manifold.

Let $(\widetilde M,\widetilde g)$ be a second manifold satisfying
the same assumptions, with an observation set $\widetilde V$
and the same fractional orders $\alpha$ and $\beta$.
We identify $\widetilde V$ with $V$ through a prescribed isometry
and use this identification for sources and observations.

\subsubsection{Recovery of the localized resolvent}

We first recover the family $L_V(z)$ defined in
\eqref{eq:localized-resolvent}.
Choose a nonnegative, nonzero function
$\eta\in C_c^\infty((0,\infty))$ and set
$f(t)=\eta(t)\phi$ for $\phi\in C_c^\infty(V)$.
By \eqref{eq:laplace-source-to-solution},
\[
    \int_0^\infty e^{-st}(L_Vf)(t)\,dt
    =
    \widehat\eta(s)L_V(-s^\beta)\phi,
    \qquad s>0.
\]
Since
\[
    \widehat\eta(s)
    =
    \int_0^\infty e^{-st}\eta(t)\,dt>0,
\]
the map $L_V$ determines $L_V(-s^\beta)\phi$ for every $s>0$.
As $-s^\beta$ ranges over $(-\infty,0)$, density of
$C_c^\infty(V)$ in $L^2(V)$ determines $L_V(z)$ as a bounded
operator for every $z<0$.

For $\phi,\psi\in C_c^\infty(V)$, define
\[
    \begin{aligned}
        m_{\phi,\psi}(z)
        &:=
        \langle L_V(z)\phi,\psi\rangle_{L^2(V)}\\
        &=
        \left\langle
            \bigl((-\Delta_g)^\alpha-z\bigr)^{-1}e_V\phi,
            e_V\psi
        \right\rangle_{L^2(M)}.
    \end{aligned}
\]
The nonnegativity and self-adjointness of $(-\Delta_g)^\alpha$
imply that $m_{\phi,\psi}$ is holomorphic on the connected domain
$\mathbb C\setminus[0,\infty)$.
Its values on $(-\infty,0)$ therefore determine it uniquely
throughout this domain.

\subsubsection{Recovery of the local spectral measure}

Define the compressed spectral measure of $(-\Delta_g)^\alpha$ by
\[
    P_V^\alpha(B)
    :=
    r_VP_{(-\Delta_g)^\alpha}(B)e_V,
    \qquad B\subset[0,\infty)\ \text{Borel}.
\]
For fixed $\phi,\psi\in C_c^\infty(V)$, let
\[
    \nu_{\phi,\psi}^\alpha(B)
    :=
    \nu_{e_V\phi,e_V\psi}^{(-\Delta_g)^\alpha}(B)
    =
    \langle P_V^\alpha(B)\phi,\psi\rangle_{L^2(V)}.
\]
If $b>0$ is not an atom ($\nu_{\phi,\psi}^\alpha(\{b\})=0$),
Stone's formula gives
\begin{equation}\label{eq:local-stone-formula}
    \nu_{\phi,\psi}^\alpha([0,b])
    =
    \lim_{\varepsilon\downarrow0}
    \frac{1}{2\pi i}
    \int_{-1}^{b}
    \bigl(
        m_{\phi,\psi}(\mu+i\varepsilon)
        -
        m_{\phi,\psi}(\mu-i\varepsilon)
    \bigr)\,d\mu.
\end{equation}
The lower endpoint lies below the spectrum of
$(-\Delta_g)^\alpha$.
Since a finite complex measure has at most countably infinite atoms, right-continuity of the cumulative measure determines
$\nu_{\phi,\psi}^\alpha([0,b])$ at all remaining $b\geq0$.
Thus $\nu_{\phi,\psi}^\alpha$ is determined on every Borel
subset of $[0,\infty)$.
Varying $\phi$ and $\psi$ and using density in $L^2(V)$
determines $P_V^\alpha$.

Let
\[
    P_V(B):=r_V P_{-\Delta_g}(B)e_V
\]
be the compressed spectral measure of $-\Delta_g$.
Since $(-\Delta_g)^\alpha$ is the fractional power of
$-\Delta_g$, the Borel functional calculus gives
\[
    P_{(-\Delta_g)^\alpha}(B)
    =
    P_{-\Delta_g}\bigl(\{\lambda\geq0:\lambda^\alpha\in B\}\bigr).
\]
Consequently,
\[
    P_V^\alpha(B)
    =
    P_V\bigl(\{\lambda\geq0:\lambda^\alpha\in B\}\bigr),
\]
and hence
\begin{equation}\label{eq:recover-local-spectral-measure}
    P_V([0,\lambda])
    =
    P_V^\alpha([0,\lambda^\alpha]),
    \qquad \lambda\geq0.
\end{equation}
This determines $P_V$ from $P_V^\alpha$.
Both measures are positive operator-valued measures; they are
not, in general, projection-valued on $L^2(V)$.

\subsubsection{Reduction to the wave equation}

For $t\geq0$, define
\[
    q_t(\lambda)
    :=
    \begin{cases}
        \displaystyle\frac{\sin(t\sqrt{\lambda})}{\sqrt{\lambda}},
            & \lambda>0,\\[6pt]
        t, & \lambda=0.
    \end{cases}
\]
The bound $|q_t(\lambda)|\leq t$ for $\lambda\geq0$ ensures
that the localized wave propagator
\begin{equation}\label{eq:localized-wave-propagator}
    K_V(t)
    :=
    r_Vq_t(-\Delta_g)e_V
    =
    \int_{[0,\infty)}q_t(\lambda)\,dP_V(\lambda)
\end{equation}
is a bounded operator on $L^2(V)$.
The last integral is understood in the weak operator sense.
Thus the recovered measure $P_V$ determines $K_V(t)$ for
every $t\geq0$.

For $F\in C_c^\infty((0,\infty)\times V)$, let $v^F$ solve
\[
    \begin{cases}
        (\partial_t^2-\Delta_g)v^F=e_VF
            & \text{in }(0,\infty)\times M,\\
        v^F|_{t=0}=0,\qquad
        \partial_tv^F|_{t=0}=0.
    \end{cases}
\]
Duhamel's formula yields
\[
    v^F(t)
    =
    \int_0^t q_{t-\tau}(-\Delta_g)e_VF(\tau)\,d\tau.
\]
Restricting to $V$, we obtain
\[
    r_Vv^F(t)
    =
    \int_0^t K_V(t-\tau)F(\tau)\,d\tau.
\]
It follows that $L_V$ determines the local source-to-solution
map for the wave equation for all positive times.
The manifold determination theorem in \cite{Correlation}
then determines $(M,g)$ up to a Riemannian isometry extending
the prescribed identification of $V$.

\bibliographystyle{unsrt}
\bibliography{references}
\end{document}